\documentclass[11pt]{article}

\usepackage{latexsym}
\usepackage{amssymb}
\usepackage{amsmath}
\usepackage{layout}
\usepackage{enumerate}
\usepackage{amsmath,amsthm,amssymb,amscd,epsfig,enumerate,graphicx}
\usepackage{makeidx}
\usepackage[all]{xy}
\usepackage{latexsym}
\usepackage{amssymb}  % Para usar los s\'{\i}mbolos de AMS
\usepackage{graphicx}
\usepackage{amsthm}
\usepackage{epsfig}

\newtheorem{tma}{Theorem}[section]
\newtheorem{lema}{Lemma}[section]
\newtheorem{prop}{Proposition}[section]
\newtheorem{cor}{Corollary}[section]

\newtheorem{defi}{Definition}[section]

\newtheorem{rem}{Remark}
\numberwithin{equation}{section}

\newcommand{\N}{\mathbb{N}}
\newcommand{\R}{\mathbb{R}}

\newcommand{\C}{\mathbb{C}}

\newcommand{\law}{\mathcal{L}}

\newcommand{\Var}{\textrm{Var}}

\newcommand{\w}{\omega}
\newcommand{\corr}{\textrm{Corr}}

\begin{document}

\begin{center}
  {\Large\bf
  Continuity in law with respect to the Hurst parameter of the
  local time of the fractional Brownian motion
  }\\*[7pt]
{ Maria Jolis}\hspace{.5 cm}{ No\`{e}lia Viles}\\
{ Departament de Matem\`atiques,  Universitat Aut\`onoma de
Barcelona\\ 08193-Bellaterra, Barcelona, Spain.\\
E-mail addresses:  {\it mjolis@mat.uab.es}}\hspace{.5 cm}{\it nviles@mat.uab.es}\\
\end{center}

\begin{abstract}
\noindent  We give
 a result of stability in law of the local time of the fractional Brownian motion  with
respect to small perturbations of the Hurst parameter. Concretely,
we prove that the law (in the space of continuous functions) of the
local time   of the fractional Brownian  motion with Hurst parameter
$H$ converges weakly  to that of the  local time of $B^{H_0}$, when
$H$ tends to $H_0$.
\end{abstract}


\section{Introduction}
In this work we will prove the continuity in law with respect to
the Hurst parameter of the family of laws of the local times of
the fractional Brownian motions.

Recall that a fractional Brownian motion $B^{H}=\{B^{H}_t,\,\,t\ge
0\}$ with Hurst parameter $H$ is a centered Gaussian process whose
covariance function is given by
$$R_H(s,t)=E[B^H_s\,B_t^H]=\frac12 \,(\,t^{2H}+s^{2H}-|t-s|^{2H}\,).$$

It is an easy exercise to see that, for any $T>0$, the family of
laws of the fractional Brownian motions $\{B^H,\,H\in (0,1)\}$
converges in law in $\mathcal C([0,T])$ to that of $B^{H_0}$, as
$H$ tends to $H_0\in (0,1)$. In fact, from the equality
$$E(B_t^H-B_s^H)^{2m}=\frac{(2m)!}{2^m\,m!}\,|t-s|^{2Hm}$$
one can see, by using Billingsley criterion (see Theorem 12.3 of
\cite{billingsley} ). that, for any $H_1\in(0,1)$, the family of
laws of $\,\{B^H,\,\,H\in[H_1,\,1)\}\,$ is tight. On the other
hand, it is clear that the covariance function $R_H(s,t)$
converges to $R_{H_0}(s,t)$, for any $s,\,t\in [0,T]$ as $H\to
H_0$, and from this we obtain the convergence of the finite
dimensional distributions.

It seems interesting to study similar results of convergence for
some functionals of the fractional Brownian motion. This kind of
results justifies the use of $B^{\hat{H}}$ as a model in
applications, when the actual value of the Hurst parameter is
unknown and $\hat{H}$ is some estimation of it.

Concretely, in this paper we consider, for any $T,\,D>0$, the
family of laws in $\mathcal C([-D,D]\times [0,T])$ of the family
$\{L^H,\,\,H\in (0,1)\}$, where $L^H=\{L^{H,t}_x,\, t\ge
0,\,x\in\R\}$ is the local time of $B^H$. We  prove, that this
family of laws converges weakly, as $H$ tends to $H_0\in(0,1)$, to
the law of $L^{H_0}$.

We point out that the existence of a  continuous version of the
local time for the fractional Brownian motion was first proved by
Berman (see \cite{berman3}). In fact, for the proof of tightness
we will mainly use the techniques based on Fourier transforms
developed by this author, jointly with a study of the correlation
of the increments of $B^H$, when $H$ belongs to a neighborhood of
$H_0$.

We have organized the paper as follows.  In the first section of
preliminaries we give the main definitions and results on the
existence and continuity of the local time for Gaussian processes
with stationary increments. In Section 2 we prove the tightness of
the laws of $\{L^H\}$ with $H$ belonging to a certain neighborhood
of $H_0$. In Section 3 we prove some general results on the
convergence in law of local times and obtain as a corollary the
desired convergence of the laws of the local times of the
fractional Brownian motions. Finally, we have added an appendix
with the proof  of a result used  in order to assure the existence
of the local time as a limit in quadratic mean.
\section{Preliminaries}

\begin{defi}\label{defmesuraocup}

Given a measurable stochastic process $X= \{X_t,\,\, t\in
[0,T]\}$, the \textsl{occupation measure} of  $X$ until the
instant $t\in [0,T] $ is defined as the following finite measure
\begin{equation*}
  \mu^t(A)=\int_0^t 1_{\{X_s\in A\}}ds,\qquad  \forall A\in\mathcal{B}(\R).
\end{equation*}
\end{defi}

\begin{defi}
A  {\it local time} of the process $X$ will be a two-parameter
process $L=\{L_{x}^t,\; x\in \R,\; t\in[0,T]\}$ such that
 for any $t\in[0,T]$ $\omega$-a.s. $L_{\cdot}^t(\omega)$ is a version of
 the density (with respect to the Lebesgue measure) of the occupation measure
$\mu^t$, in the case in which this density exists.
\end{defi}

\begin{rem}
Notice  that $L_{\cdot}^0$ can be taken identically equal to $0$
and that the existence of a density for $\mu^T$ implies the
existence of a density for $\mu^t$ for any $t\in [0,T]$.
\end{rem}

One of the main properties of the occupation measure is given by
the following result, known as {\it occupation formula} (see for
instance  \cite{berman1}).

\begin{prop}\label{prop12}
Let $X=\{X_t,\;\, t\in [0,T]\}$ be a measurable process. If $g: \R
\rightarrow \C$ is a Borel measurable function, then
\begin{equation*} \int_{\R} g(u) d\mu^t(u)=\int_0^t g(X_s)ds,
\end{equation*}
and both integrals are well defined or not at the same time.
\end{prop}

We will use the following version of Plancherel's theorem.

\begin{tma}\label{plancherel}  If the Fourier transform $\phi$ of a measure $\mu$ belongs to
$L^2(\R)$, then $\mu$ is absolutely continuous and its density is
also square integrable. Moreover the density,  $f$, is the limit
in $L^2(\R)$ of
\begin{equation}\label{eq:eq6}
f_N(x)=\frac{1}{2\pi}\int_{-N}^N e^{-iux}\phi(u)du,
\end{equation}
when $N\to \infty$.
\end{tma}

\begin{rem}
We denote by $\,\phi^t$ the Fourier transform of the occupation
measure $\mu^t$. Then, by using the occupation formula
(Proposition \ref{prop12}), $\phi^t\in L^2(\R)$ almost surely,  if
and only if
$$\int_{\R}\Big|\int_0^t e^{iuX_s}ds\Big|^2 du<\infty\;\; a.s.,$$
or equivalently, if
$$
\int_\R\left(\int_0^t\int_0^t
e^{iu(X_s-X_r)}drds\right)du<\infty\;\;a.s.
$$
\end{rem}

\bigskip

\noindent The fact that the limit in (\ref{eq:eq6}) is in $L^2(\R)$
has the inconvenience that, in this case, the density $f(x)$ is
defined  for all $x$-a.e. By this reason we will need some results
that allow to ensure the existence of
 $L=\{L_x^t, \; x\in \R,\;t\in
[0,T]\}$ as a stochastic process. The following result is given in
Theorem 4.1 of \cite{berman1}:

\begin{tma}\label{thmabermanl2}
Suppose that
\begin{equation*}
\int_{\R}\int_{\R}\int_0^T\int_0^T|E[e^{i u X_s+i v X_r}]|ds dr du
dv<+\infty.
\end{equation*}
Define
$$\psi_N(x,t,\w)=\frac{1}{2\pi}\int_{-N}^N
e^{-iux}\int_0^te^{iu X_r(\omega)}dr du.$$ Then, for any
$(x,t)\in\R\times [0,T]$ there exists a random variable  $L_x^t$
such that
\begin{equation*}
\lim_{N\to \infty}\sup_{(x,t)\in
\R\times[0,T]}E|\psi_N(x,t)-L_x^t|^2=0.
\end{equation*}

\end{tma}

The following theorem allows us to obtain the local time as a
limit in quadratic mean. Its proof, that uses  Theorem
\ref{thmabermanl2} and Plancherel's theorem, is given in the
Appendix.

\begin{tma}\label{extlocal}
Let $X=\{X_t, \; t\in [0,T]\}$ be a measurable stochastic process
 verifying the following conditions:
\begin{enumerate}[(i)]
\item $$\int_\R\int_\R\int_0^T\int_0^T |E[e^{iuX_s+ivX_r}]|ds
drdudv<\infty.$$ \item For each $t\in [0,T]$,
$$\int_\R\int_0^t\int_0^t
e^{iu(X_s-X_r)}drdsdu<\infty,\quad\quad \mbox{a.s\,}.$$
\end{enumerate}
Consider $L=\{L_x^{t}, \; (x,t)\in\R\times[0,T]\}$ the process
defined for any $(x,t)$ as the random variable $L_x^t$ appearing
in Theorem  \ref{thmabermanl2}. Then, this process $L$ is a local
time of $X$.
\end{tma}

The next lemma (see Example 3.2 of \cite{berman1}) gives a
sufficient condition in order that a Gaussian process with
stationary increments satisfies hypothesis (ii) of Theorem
\ref{extlocal}.

\begin{lema}\label{mesocupl2}
Let $X=\{X_t, \;t\in [0,T]\}$ be a measurable centered Gaussian
process  null at  $0$ with stationary increments. Denote  by
$\sigma^2(t)$ the  variance function of $X$. If
\begin{equation}\label{lema323232}
\int_0^T\int_0^t \sigma(t-s)^{-1}dsdt<\infty,
\end{equation}
 then, for any  $t\in[0,T]$,
\begin{equation*}\label{transfL2}
E\left(\int_\R\int_0^t\int_0^t
e^{iu(X_s-X_r)}drdsdu\right)<\infty.
\end{equation*}

\end{lema}

We will also need the following lemma.

\begin{lema}\label{desvar} (see Lema 8.1 of
\cite{berman5}) Let $Y_1,\dots, Y_m$ be non-constant square
integrable random variables. The following inequality is satisfied
for all  $v_1,\dots,v_m\in\R$:
\begin{equation*}
Var\left(\sum_{i=1}^m v_iY_i \right)\geqslant \frac{\det \Gamma
}{\prod_{i=1}^m\Gamma_{i i} }\frac{1}{m}\sum_{i=1}^m
v_i^2\Gamma_{i i },
\end{equation*}
where  $\Gamma$ is the covariance matrix of $\,Y_1,\dots, Y_m$.
\end{lema}

\bigskip

\noindent We will also use the equality given in the following
lemma.

\begin{lema}\label{int}
For any  $a>0$ and $0<\alpha<2$,
$$\int_{\R} |x|^{\alpha} e^{-ax^2}dx=a^{-\frac{(\alpha+1)}{2}}\Gamma\left(\frac{\alpha+1}{2}\right).$$
\end{lema}

\bigskip

\noindent Given a function $F$, defined on  $\R^2$, and
$(s,t),\;(s',t')\in \R^2$ such that $s\le s'$ and $t\le t'$, we
will denote by $\Delta_{s,t}F(s',t')$ the increment of $F$ over
the rectangle $((s,t),(s',t')]$, that is,
$$\Delta_{s,t}F(s',t')=F(s',t')-F(s',t)-F(s,t')+F(s,t).$$

The next theorem is the main result of this section where
sufficient conditions are given for a Gaussian process with
stationary increments  to have a local time possessing a
continuous version. This theorem is an adaptation of Theorem 8.1
of  \cite{berman5}. We will give its proof because we will need a
precise evaluation of the constants appearing in it.

\begin{tma}\label{tmamomcom}
Let $X=\{X_t, t\in [0,T]\}$ be a centered Gaussian measurable
process  null at
 $0$ with  stationary increments such that its variance function  $\sigma^2(t)$ is bounded by a constant $C_{\sigma}$.
  Suppose that
\begin{enumerate}[(i)]
\item There exists $m_0$, even natural number, such that for any $m\geqslant m_0$ even, the determinant
of the covariances of the normalized increments
$$\frac{X_{t_j}-X_{t_{j-1}}}{\sigma(t_j-t_{j-1})},\qquad\qquad \qquad j=1,\dots, m,$$
is bounded from below by a constant $A_m>0$ on the set
\\$\{(t_1,\dots,t_m)\in [0,T]^m: 0=t_0<t_1<\dots<t_m<T\}$.
\item  There exist $\delta>0$ and $\alpha>0$ such that
\begin{equation}\label{condiithma105}
\sup_{t\in [0,T]}\int_t^{t+h}
[\sigma(s)]^{-(1+2\delta)}ds\leqslant C_{\alpha,\delta}\,
h^{\alpha}.
\end{equation}
\end{enumerate} Then,
\begin{enumerate}[a)]
\item For each $(x,t)$, there exists the local time  $L_x^t$ as a
limit (uniform in $(x,t)$) in quadratic mean.
\item For any even $m\geqslant
m_0$, there exists a positive constant  $C_1$ depending on
 $m$, $A_m$, $\alpha$ and $\delta$ such that
\begin{equation*}
E|\Delta_{0,t}L(0,t+h)|^m\leqslant C_1 |h|^{m\alpha}.
\end{equation*}
We can take  $C_1=C_m A_m^{-m/2}(C_{\alpha,\delta})^m$, with $C_m$
only depending on $m$.
\item If $\,m\geqslant m_0$ and even, there exists a positive constant
$C_2$ that depends on $m$, $A_m$, $\delta$, $\alpha$ and $\sigma$
such that
\begin{equation*}
E|\Delta_{x,t}L(x+k,t+h)|^m\leq C_2|h|^{m\alpha}|k|^{m\delta}.
\end{equation*}
We can take
\begin{equation*}
C_2=C_m\max(1,A_m^{-m/2})(\max(1,C_\sigma^{2\delta}))^mC_{\alpha,\delta}^m,
\end{equation*}
with $C_m$ only depending on $m$.
\end{enumerate}

\noindent As a consequence of $b)$ and  $c)$, by using
  Kolmogorov-Chentsov's criterion, we obtain the existence of a
 version of the local time of $X$, $L=\{L_x^{t}, \;
(x,t)\in\R\times [0,T]\}$ that is jointly continuous in $(x,t)$.
\end{tma}

\begin{proof}
First of all, we check that the hypotheses of this theorem imply
those of Theorem  \ref{extlocal}. This will give  $(a)$.

Indeed, it is easy to see that condition  $(ii)$  implies
inequality (\ref{lema323232}) of Lemma \ref{mesocupl2} and, as a
consequence, hypothesis  $(ii)$ of Therem \ref{extlocal} is
satisfied.

On the other hand, it is not difficult to check that  for  $m\ge
2$ even,
$$
\aligned \int_0^{T}\int_0^{T}&\int_{-\infty}^\infty
\int_{-\infty}^\infty E[e^{iuX_s+ivX_t}]\,ds\,dt\,du\,dv \\
& \le
\left(\int_0^{T}\!\!\!\!\cdots\!\!\!\!\!\int_0^{T}\!\int_{-\infty}^\infty
\!\! \dots \!\! \int_{-\infty}^\infty E[e^{\sum_{j=1}^m
iu_jX_{s_j}}]\prod_{j=1}^m du_j \prod_{j=1}^m ds_j\right)^{2/m}.
\endaligned$$
This last integral is finite for  $m\ge m_0$, by the  arguments
that we will give below. This provides us condition  $(i)$ of
Theorem \ref{extlocal}.

Both  $(b)$ and $(c)$ are proved in a similar way, we will only
give the proof of $(c)$ that is the more complicated one.

Taking into account  $(a)$, for $m$ even, we can express the
$m$-th moment of the $2$-dimensional increment of the local time
of $X$ as
$$
\aligned E&|\Delta_{x,t}L(x+k,t+h)|^m =E[\Delta_{x,t}L(x+k,t+h)]^m\\
& = (2\pi)^{-m} E\left(\lim_{N\to
\infty}\int_{-N}^N\!\!\cdots\int_{-N}^N\int_t^{t+h}\!\!\!\cdots\int_t^{t+h}
 \prod_{j=1}^m
(e^{-iu_j(x+k)}-e^{-iu_jx})\right.\\& \left.\times\prod_{j=1}^m
e^{iu_jX_{s_j}}\prod_{j=1}^m ds_j \prod_{j=1}^m du_j \right).
\endaligned
$$
It can be checked that the above expression  can be bounded by
\begin{equation*}
(2\pi)^{-m}\!\!\!\!\int_t^{t+h}\!\!\!\!\!\!\cdots\!\!\!\int_t^{t+h}\!\!\!
\int_{-\infty}^{\infty}\!\!\!\cdots\!\!\int_{-\infty}^{\infty}
\!\!\prod_{j=1}^m
\left(|e^{-iu_j(x+k)}\!\!\!-\!\!\!e^{-iu_jx}|\right)\!\!E[e^{\sum_{j=1}^m\!
iu_jX_{s_j}}]\prod_{j=1}^mdu_j\prod_{j=1}^m ds_j.
\end{equation*}
Using that  $|e^{ix}-e^{iy}|\le 2\,|x-y|^{\delta}$, for any
$x,\,y\in\R$ and all  $\delta\in (0,1)$, we have that
\begin{equation*}
\aligned
&\int_t^{t+h}\!\!\cdots\!\!\!\int_t^{t+h}\int_{-\infty}^\infty
\!\! \dots \!\! \int_{-\infty}^\infty \prod_{j=1}^m
|e^{-iu_j(x+k)}-e^{-iu_jx}|E[e^{\sum_{j=1}^m
iu_jX_{s_j}}]\prod_{j=1}^mdu_j\prod_{j=1}^m ds_j\\&\leqslant
2^m|k|^{m\delta}
\int_t^{t+h}\!\!\cdots\!\!\!\int_t^{t+h}\int_{-\infty}^\infty \!\!
\dots \!\! \int_{-\infty}^\infty \prod_{j=1}^m |u_j|^ \delta
E[e^{\sum_{j=1}^m iu_jX_{s_j}}]\prod_{j=1}^mdu_j\prod_{j=1}^m
ds_j.
\endaligned
\end{equation*}

Given that  the integrand on the last expression is symmetric in
$s_1,\dots,s_m$, we can change the domain of integration, $[t,
t+h]^m$, by its subset
$$\{(s_1,\dots,s_m): t\le s_1<\dots<s_m\le t+h\}.$$
Making the following change of variables
\begin{equation*}\label{canvivar} \left\{
\begin{array}{ll}
u_j=v_j-v_{j+1}, & \mbox{$\forall j=1,\dots, m-1,$}\\
u_m=v_m ,&
\end{array} \right.
\end{equation*}
and defining  ${s_0}=0$, we have that
$$\sum_{j=1}^m u_jX_{s_j}=\sum_{j=1}^m v_j (X_{s_j}-X_{s_{j-1}}),
$$
and this entails that
$$E[e^{\sum_{j=1}^m i\,u_j\,X_{s_j}}]=e^{-\frac{1}{2}\Var[\sum_{j=1}^m v_j(X_{s_j}-X_{s_{j-1}})] }.$$

By Lemma  \ref{desvar}, we can majorize this expression as follows
$$E[e^{\sum_{j=1}^m iv_j(X_{s_j}-X_{s_{j-1}})}]\leqslant
\exp \left(-\frac{1}{2}\frac{R}{\prod_{j=1}^m
\sigma^2(s_j-s_{j-1})} \frac{1}{m}\sum_{j=1}^m
v_j^2\sigma^2(s_j-s_{j-1})\right),$$ where $R$ is the determinant
of the covariance matrix of the increments $X_{s_j}-X_{s_{j-1}}$
for $j=1,\ldots, m-1$.

Since the constant  $A_m$ is a lower bound of the determinant of
the correlation matrix of the increments of  $X$ that coincides
with $\frac{R}{\prod_{j=1}^m \sigma^2(s_j-s_{j-1})}$, we have that
\begin{equation*}
E[e^{i\sum_{j=1}^m v_j(X_{s_j}-X_{s_{j-1}})}]\leqslant
e^{-B_m\sum_{j=1}^m v_j^2 \sigma^2(s_j-s_{j-1})},
\end{equation*}
with $B_m=\frac{A_m}{2m}$.

On the other hand,
\begin{equation*}
\prod_{j=1}^m |u_j|^{\delta}=\left(\prod_{j=1}^{m-1}
|v_j-v_{j+1}|^{\delta}\right)|v_m|^{\delta}\leqslant
\left(\prod_{j=1}^{m-1}(|v_j|^\delta +
|v_{j+1}|^\delta)\right)|v_m|^{\delta}.
\end{equation*}
This last product equals to the sum of $2^{m-1}$ terms, each of
them containing at most  $m$
factors $|v_1|^\delta,\dots, |v_m|^\delta$ with exponents  $0$, $1$ or $2$.\\
Using this, we obtain
\begin{align*}
&E|\Delta_{x,t}L(x+k,t+h)|^m\\&\leqslant
|k|^{m\delta}\pi^{-m}m!\\&\!\!\times\!\!\!\sum_{\theta_i\in
\{0,1,2\}}\int_{\{t\leqslant s_1< \dots < s_m \le t+h\}\times
\R^m}\!\!\!\!\! |v_1^{\theta_1}\dots v_m^{\theta_m}|^\delta
e^{-B_m\sum_{j=1}^m v_j^2 \sigma^2(s_j-s_{j-1}) }  \prod_{j=1}^m
dv_j \!\!\prod_{j=1}^m ds_j.
\end{align*}
By  Fubini's theorem and Lemma \ref{int},
\begin{align}\label{eq:eq210}
&E|\Delta_{x,t}L(x+k,t+h)|^m \nonumber\\ &\leqslant
|k|^{m\delta}\pi^{-m}m!\!\!\!\!\sum_{\theta_i\in \{0,1,2\}}
{\idotsint}_{\substack{\{t\leqslant s_1< \dots< s_m \le
t+h\}}}\!\!\prod_{j=1}^m
(B_m\sigma^2(s_j-s_{j-1}))^{-\frac{(\theta_j\delta+1)}{2}}\nonumber\\&\times\prod_{j=1}^m\Gamma
\left(\frac{\theta_j \delta+1} {2}\right)ds_1\dots ds_m.
\end{align}
Taking into account that $
\max_{r\in[\frac{1}{2},\frac{3}{2}]}\Gamma(r)=\Gamma(\frac12)=\sqrt{\pi}
$, we can  bound expression (\ref{eq:eq210}) in the following way
\begin{align*}
E|\Delta_{x,t}&L(x+k,t+h)|^m \leqslant|k|^{m\delta}\pi^{-m/2}\,
m!\,\max(1,B_m^{-m/2})
\qquad\qquad\qquad\\&\times\sum_{\theta_j\in\{0,1,2\}}{\idotsint}_{\substack{\!\!\!\{t\leqslant
s_1< \dots< s_m \le t+h\}}}\prod_{j=1}^m
[\sigma(s_j-s_{j-1})]^{-(\theta_j\delta+1)}ds_1\dots ds_m.\\
&\leqslant |k|^{m\delta}\pi^{-m/2}\,
m!\,\max(1,B_m^{-m/2})\left(\max(1,C_\sigma^{2\delta})\right)^m\qquad\nonumber\\&\!\!\times\!\!
{\idotsint}_{\!\!\substack{\!\!\!\!\{t\leqslant s_1< \dots< s_m
\le t+h\}}}\prod_{j=1}^m
[\sigma(s_j-s_{j-1})]^{-(1+2\delta)}ds_1\dots ds_m.
\end{align*}

Finally, using that
\begin{align*}
{\idotsint}_{\substack{\{t\leqslant s_1< \dots< s_m \le
t+h\}}}\!\!&\prod_{j=1}^m
[\sigma(s_j-s_{j-1})]^{-(1+2\delta)}ds_1\dots ds_m \nonumber\\&
\\&\!\!\!\!\!\!\!\!\!\!\le \left(\int_t^{t+h}
\sigma(x)^{-(1+2\delta)}dx\right)\left(\int_0^h
\sigma(x)^{-(1+2\delta)}dx\right)^{m-1},
\end{align*}
 denoting by    $C_m$ the product of all the constants depending only on  $m$ and using
 (\ref{condiithma105}), we have that
\begin{align*}
E|\Delta_{x,t}L(x+k,t+h)|^m\leqslant
&C_m\max(1,A_m^{-m/2})\left(\max(1,C_\sigma^{2\delta})\right)^m(C_{\alpha,\delta})^m
h^{m\alpha}|k|^{m\delta}.
\end{align*}
This ends the proof.
\end{proof}
\section{Existence and continuity of the local time of the fractional Brownian motions. Tightness of
the family of their laws}

Recall that  the fractional Brownian motion of Hurst parameter
$H\in (0,1)$, denoted by $B^H$, is a centered  Gaussian process
with stationary increments, and  taking a continuous version we
ensure that it is also a measurable process.

Along this section we will check that the fractional Brownian
motion satisfies the conditions of Theorem \ref{tmamomcom}.  And
we will see that all the constants appearing in these conditions
can be taken independent of the parameter $H$, at least in a
neighborhood of $H_0$, for any $H_0\in(0,1)$.

First of all, notice that the variance function of the fractional
Brownian motion of parameter $H$, $\sigma_H^2(t)=t^{2H}$ is
bounded by  $C_{T}=\max(1,T^{2})$ for $t\in [0,T]$. Then, the
constant $C_\sigma$ appearing in Theorem \ref{tmamomcom} equals to
$C_T$.

We state in the next lemma, whose proof is a simple computation,
that the variance function $\sigma_H^2$ also satisfies condition
(\ref{condiithma105}) of Theorem \ref{tmamomcom}.
\begin{lema}\label{cotintsigma}
Let $H\in(H_0-\eta,H_0+\eta)\subset (0,1)$. Then, for any  $\delta
>0$ satisfying $(H_0+\eta)(1+2\delta)<1$,
\begin{equation*}
\int_t^{t+h} [\sigma(s)]^{-(1+2\delta)}ds\leqslant
C_{T,H_0,\eta,\delta}\; h^{1-(H_0+\eta)(1+2\delta)},
\end{equation*}
where
$$C_{T,H_0,\eta,\delta}=\frac{\max(1,T^{2\eta(1+2\delta)})}{(1-(H_0+\eta)(1+2\delta))T^{1-(H_0+\eta)(1+2\delta)}}.$$
\end{lema}

Now, we will prove that for any  $m\geqslant 2$, the determinant
of the covariance matrix of the normalized increments of the
process is bounded from below by a positive constant $A_m^H$ over
the set $\{(t_1,\dots,t_m)\in[0,T]^m:0=t_0<t_1<\dots<t_m<T\}$.
Moreover, we will show that this constant can be taken
independently of $H$, at least in a neighborhood of a point
$H_0\in (0,1)$.

We will need to distinguish the cases $H_0<\frac12$ and $H_0\ge
\frac12$. For $H_0\in(0,\frac{1}{2})$, we will use the following
lemma.

\begin{lema}\label{cotdet}(see \cite{berman3})
Let $X=\{X_t, t\in[0,T]\}$ be a Gaussian process with stationary
increments and concave variance function. Let $0\leqslant
t_0<t_1<\!\!\dots\!<\!t_n$. Then the following inequality is
satisfied
\begin{equation*}
\frac{\det(E[(X_{t_i}-X_{t_{i-1}})(X_{t_j}-X_{t_{j-1}})])_{1\leqslant
i,j\leqslant n}}{\prod_{j=1}^n \sigma^2(t_i-t_{i-1}) }\geqslant
2^{-3n},
\end{equation*} for $n$ positive integer.
\end{lema}

Since the variance function of the fractional Brownian motion with
parameter $H\in (0,\frac12]$ is concave, we can apply this lemma
and obtain that the constant $A_m^H$ equals to $2^{-3m}$ (notice
that, in fact, for $H=\frac12\,$ we can take $A_m^H=1$ ).

When $H\in(\frac{1}{2},1)$, the above result does not apply. In
fact, we will need to analyze the behaviour of the elements of the
covariance matrix of the normalized increments of $B^H$ in a small
neighborhood of $H_0\in [\frac12, 1)$.

The correlation  of two disjoint increments of $B^H$ is given by
\begin{equation*} \corr(B^H_t-B^H_s,
B^H_v-B^H_u)=\frac{1}{2}\,\frac{(u-t)^{2H}-(u-s)^{2H}-(v-t)^{2H}+(v-s)^{2H}}{(t-s)^H(v-u)^H},
\end{equation*}
with  $0\leqslant s<t\leqslant u<v\leqslant T$.  We know that, for
$H>\frac12$, this correlation is always positive.

If we write
\begin{align*}
v-u&=\gamma(t-s),\\
u-t&=\beta(t-s),
\end{align*}
we will have
\begin{equation*}
\corr(B^H_t-B^H_s,
B^H_v-B^H_u)=\frac{1}{2}\,\frac{\beta^{2H}-(1+\beta)^{2H}-(\beta+\gamma)^{2H}+(1+\beta+\gamma)^{2H}}{\gamma^H}.
\end{equation*}
If we consider consecutive increments $B^H_t-B^H_s$ and
$B^H_w-B^H_t$, therefore $u=t$ (this entails that $\beta=0$) and
we have
\begin{equation*}
\corr(B^H_t-B^H_s,B^H_w-B^H_t)=\frac12\,\frac{(1+\gamma)^{2H}-\gamma^{2H}-1}{\gamma^H}.
\end{equation*}

\begin{lema}\label{lemainconsec}
Let $H\in(\frac{1}{2},1)$. For any  $\gamma, \beta$ positive real
numbers, the following inequality is satisfied
\begin{equation*}
\frac{\beta^{2H}-(1+\beta)^{2H}-(\beta+\gamma)^{2H}+(1+\beta+\gamma)^{2H}}{\gamma^H}\leqslant
\frac{(1+\gamma)^{2H}-\gamma^{2H}-1}{\gamma^H}.
\end{equation*}
\end{lema}
\begin{proof}
The proof is a simple argument of convexity.
\end{proof}
In the following proposition we will prove that condition $(i)$ of
Theorem  \ref{tmamomcom} is satisfied and that the constant
$A_m^H$ can be taken independently of $H$ for $H$ belonging to a
neighborhood of $H_0\ge \frac12$.
\begin{prop}\label{detposhm12}
Let $H_0\in[\frac{1}{2},1)$. For any $m\ge 2$ there exist
$\eta>0$ and a constant  $A_{m,\eta,H_0}>0$ depending only on
 $m$, $\eta$ and $H_0$, such that the determinant of the correlation matrix of
 the increments $B^H_{t_j}-B^H_{t_{j-1}}$,
$j=1,\;\dots,\;m$,  on the set $\{(t_1,\!\dots\!,t_m)\in
[0,T]^m:\; 0=t_0\!<\!\dots\!<t_m<T\}$, \linebreak is greater or
equal than  $A_{m,\eta,H_0}$, for any $H\in (H_0-\eta,H_0+\eta)$.
\end{prop}
\begin{proof}
The existence, for any $H\in (\frac12, 1)$, of a constant
$A_{m,H}>0$ such that it is a lower bound for
 the determinant of the correlation matrix of the increments of $B^H$ is proved in Theorem
6.2 of \cite{berman3}.

 In order to prove the proposition, we will show that in a small
 enough
 neighborhood of $H_0$ the determinant of the correlation matrix
 corresponding to $B^H$ is near to the corresponding to
$B^{H_0}$.

Notice that if $H_0=\frac{1}{2}$,  Lemma \ref{cotdet} says us that
the determinant of the correlation matrix of the increments of
$B^H$ is bounded from below by $2^{-3m}$  for
$H\in(H_0-\eta,H_0]$, for any $\eta<\frac12$.  The arguments that
we will use from now on
 will provide us a neighborhood of the form
$(H_0,H_0+\eta)$ in which the uniform lower bound of the
determinant also exists.
% D'aquesta manera usant el resultat Berman que hem citat anteriorment obtindrem una constant
%$A_{H_0}$ propera a $A_H$ i tamb\'{e} positiva.

Since the determinant of a matrix is a sum of products of its
components, and taking into account that, for any $H\in(0,1)$,
\begin{equation*}\label{fitcorrel1}
\sup_{0\leqslant s< t\leqslant u<v\leqslant
T}|\corr(B^{H}_t-B^{H}_s,B^{H}_v-B^{H}_u)|\leqslant 1\,,
\end{equation*}
the proof will be concluded if we see that for any $\varepsilon
>0$ there exists  $\rho>0$ such that for any $H\in (H_0-\rho,
H_0+\rho)$ we have that
\begin{equation}\label{difcorrel}
\sup_{0\le s< t\leqslant u<v\leqslant
T}|\corr(B^H_t-B^H_s,B^H_v-B^H_u)-\corr(B^{H_0}_t-B^{H_0}_s,B^{H_0}_v-B^{H_0}_u)|<\varepsilon.
\end{equation}

% Haurem vist que el determinant de la matriu de correlaci\'{o} corresponent a X^H \'{e}s proper al del de la matriu corresponent a X^{H_0}
%quan H \'{e}s propera a H_0
As we have seen above, we can express the correlation between two
disjoint increments $B^H_t-B^H_s$ and $B^H_v-B^H_u$ in terms of the
parameters  $\beta$ and $\gamma$, with
\begin{align*}
v-u&=\gamma(t-s),\\
u-t&=\beta(t-s),
\end{align*}
and  $0\leqslant s<t\leqslant u<v\leqslant T$. So,
(\ref{difcorrel}) is equivalent to
\begin{align}\label{difcorrel2}
\sup_{\substack{0<\gamma< +\infty\\0\le \beta<+\infty
}}&\Big|\frac{\beta^{2H}-(1+\beta)^{2H}-(\beta+\gamma)^{2H}+(1+\beta+\gamma)^{2H}}{\gamma^H}\nonumber
\\&-\frac{\beta^{2H_0}
-(1\!+\!\beta)^{2H_0}-(\beta+\gamma)^{2H_0}+(1+\beta+\gamma)^{2H_0}}{\gamma^{H_0}}\Big|<\varepsilon,
\end{align}
for $|H-H_0|$ small enough.

Now, we will show (\ref{difcorrel2}). Taking into account the
different possible values of the parameters $\gamma$ and $\beta$,
we will prove the following assertions.
\begin{enumerate}[(i)]
\item For any $\varepsilon>0$ and any $0<\delta<M$,
$L>0$ there exists $\rho_1$ (depending on $\varepsilon$, $M$,
$\delta$  and $L$) such that, for $|H-H_0|<\rho_1$, we have
\begin{align}\label{fcorr31}
\sup_{\substack{\delta\leqslant \gamma\leqslant M\\0\leqslant
\beta\leqslant
L}}\Big|&\frac{\beta^{2H}-(1+\beta)^{2H}-(\beta+\gamma)^{2H}+(1+\beta+\gamma)^{2H}}{\gamma^H}\nonumber\\&-\frac{\beta^{2H_0}
-(1+\beta)^{2H_0}-(\beta+\gamma)^{2H_0}+(1+\beta+\gamma)^{2H_0}}{\gamma^{H_0}}\Big|<\varepsilon.
\end{align}
This is a consequence of the uniform continuity of the function
$\!f(x,y)\!=x^y\!$ on any compact that does not contain the point
$(0,0)$.

\item Given $\varepsilon>0$, there exists $\delta>0$ and $\rho_2>0$
such that  for any  $|H-H_0|<\rho_2$, we have
\begin{equation}\label{fcorrel322}
\sup_{\substack{0<\gamma<\delta\\ 0\leqslant
\beta<+\infty}}\frac{\beta^{2H}-(1+\beta)^{2H}-(\beta+\gamma)^{2H}
+(1+\beta+\gamma)^{2H}}{\gamma^H}<\varepsilon.
\end{equation}

Notice that (\ref{fcorrel322}) gives  (\ref{difcorrel2}) when we
take the supremum on $0<\gamma<\delta$, $0\leqslant \beta
<\infty$, because (\ref{fcorrel322}) is also valid for $H=H_0$.

From Lemma \ref{lemainconsec} we know that
\begin{equation}\label{fcorr32}
\frac{\beta^{2H}-(1+\beta)^{2H}-
(\beta+\gamma)^{2H}+(1+\beta+\gamma)^{2H}}{\gamma^H}\leqslant\frac{(1+\gamma)^{2H}
-\gamma^{2H}-1}{\gamma^H},
\end{equation}
Moreover,  the two numerators of the above expressions are
positive.

By the Mean Value theorem, we have that
\begin{equation*}
(1+\gamma)^{2H}-\gamma^{2H}-1=2H((1+\xi)^{2H-1}-\xi^{2H-1})\gamma,
\qquad\qquad \xi\in(0,\gamma).
\end{equation*}
Using this inequality and taking into account that $0<2H-1<1$, we
can bound the right-hand side of  the inequality (\ref{fcorr32})
in the following way
\begin{equation}\label{cotacota}
\frac{(1+\gamma)^{2H}-\gamma^{2H}-1}{\gamma^H} \leqslant
\,2\,(H_0+\rho_2)\,\delta^{1-H_0-\rho_2}\,<\,\varepsilon,
\end{equation}
by taking  $\,\rho_2$ verifying $1-H_0-\rho_2>0$ and
$$0<\delta<\left(\frac{\varepsilon}{2(H_0+\rho_2)}\right)^{\frac{1}{1-H_0-\rho_2}}.$$
This gives (\ref{fcorrel322}).
\item Given $\varepsilon>0$, there exists $M>0$ and $\rho_3>0$ such that for
$|H-H_0|<\rho_3$, we have
\begin{equation}\label{fraccorr4}
\sup_{\substack{\gamma>M\\0\leqslant\beta<+\infty}}
\frac{\beta^{2H}-(1+\beta)^{2H}-(\beta+\gamma)^{2H}+(1+\beta+\gamma)^{2H}}{\gamma^H}<\varepsilon.
\end{equation}
Indeed, using again Lemma \ref{lemainconsec},
\begin{align*}\label{fraccorr444}
\sup_{\substack{\gamma>M\\0\leqslant\beta<+\infty}}&\frac{\beta^{2H}-(1+\beta)^{2H}
-(\beta+\gamma)^{2H}+(1+\beta+\gamma)^{2H}}{\gamma^H}\\\leqslant&
\sup_{\gamma>M}\frac{(1+\gamma)^{2H}-\gamma^{2H}-1}{\gamma^H}\\
=&\sup_{0<y<\frac1{M}}\frac{\left(1+\frac{1}{y}\right)^{2H}-\left(\frac{1}{y}\right)^{2H}-1}{\frac{1}{y^H}}
= \sup_{0<y<\frac1{M}}\frac{(1+y)^{2H}-(1+y^{2H})}{y^H}.
\end{align*}
So,  (\ref{fraccorr4}) is a consequence of (\ref{cotacota}).

\item Finally, given $\varepsilon>0$ and $0<\delta<M$, there exists
$L>0$ and $\rho_4>0$ such that for $|H-H_0|<\rho_4$, we have
\begin{equation*}
\sup_{\substack{\delta_\leqslant \gamma\leqslant M\\\beta>L}}
\frac{\beta^{2H}-(1+\beta)^{2H}-(\beta+\gamma)^{2H}+(1+\beta+\gamma)^{2H}}{\gamma^H}<\varepsilon.
\end{equation*}
Indeed, on one hand $\gamma$ belongs to the compact, $[\delta,
M]$, far away from  $0$. So, it suffices to study the numerator of
the above expression.  Applying twice the Mean Value Theorem we
obtain
\begin{align*}
0\leqslant
(1+\beta+\gamma)^{2H}-(\beta+\gamma)^{2H}-&((1+\beta)^{2H}-\beta^{2H})=2H\!(\xi^{2H\!-\!1}-\eta^{2H-1})\\&=
\!2H\!(2H\!-\!1)\upsilon^{2H-2}(\xi-\!\eta),
\end{align*}
where $\xi\in(\beta+\gamma, \beta+\gamma+1)$,
$\eta\in(\beta,\beta+1)$, $\upsilon\in\langle\eta,\xi\rangle$,
from which we deduce that $\xi\ge\eta$.

Taking into account that $\xi\in(\beta+\gamma, \beta+\gamma+1)$,
$\eta\in(\beta,\beta+1)$ and $\delta\leqslant\gamma\leqslant M$,
we have that $0\le\xi-\eta<M+1.$

On the other hand, since $L<\upsilon<\beta+M+1$ and $2H-2<0$, we
obtain
\begin{equation*}
(1+\beta+\gamma)^{2H}-(\beta+\gamma)^{2H}-((1+\beta)^{2H}-\beta^{2H})\leqslant
2H(2H-1)\left(\frac{1}{L}\right)^{2-2H}(M+1).
\end{equation*}
Finally, since we can take  $L>1$, we obtain
\begin{align*}
(1+\!\beta+\!\gamma)^{2H}-(\beta+\gamma)^{2H}-&((1+\beta)^{2H}-\beta^{2H})\leqslant
2(H_0+\rho_4)(2(H_0+\rho_4)\!-\!1)\\&\times\left(\frac{1}{L}\right)^{2-2(H_0+\rho_4)}(M+1)<\varepsilon,
\end{align*}
if $|H-H_0|<\rho_4$, with  $\rho_4>0$ satisfying
$2-2(H_0+\rho_4)>0$ (or equivalently, $\,0<\rho_4<1-H_0$) and $L$
big enough.
\end{enumerate}

This finishes the proof of (\ref{difcorrel2}).

\end{proof}

As a consequence of the previous results of this section and
Theorem \ref{tmamomcom}, we can state the following proposition.

\begin{prop}\label{ajustprev}
Let $H_0\in (0,1)$. Then, there exists $\eta>0$ such that the
family  $\{B^H,\; H\in(H_0- \eta, H_0+\eta)\}$ satisfies
\begin{enumerate}[(i)]
\item For any  $(x,t)$ and each $H\in (H_0-\eta,H_0+\eta)$,
there exists the local time $L_x^{t,H}$ as a limit (uniform in
$(x,t)$) in quadratic mean.
\item There exist  positive constants $C_1$ (depending on $m$,
$\eta$ and $H_0$) and $\alpha$ (depending on $\eta$ and $H_0$),
such that
\begin{equation*}
E|\Delta_{0,t}L^H(0,t+h)|^m\leqslant C_1 |h|^{m\alpha}.
\end{equation*}
\item For all $m$ even and $m\geqslant 2$, there exist
$\delta>0$ and $\alpha>0$ depending on $H_0$ and $\eta$ and also
$C_2>0$ depending on $m$, $H_0$ and $\eta$ such that
\begin{equation*}
E|\Delta_{x,t} L^H(x+k, t+h)|^m\leqslant
C_2|h|^{m\alpha}|k|^{m\delta}.
\end{equation*}
\end{enumerate}
As a consequence of (ii) and (iii), by using the tightness
criterion of \cite{yor}, we obtain
 the tightness of the laws of
$\{L^H,\,H\in(H_0-\eta,\,H_0+\eta)\}$ in $\mathcal C([-D,D]\times
[0,T])$, for any $D>0$ and $T>0$.
\end{prop}

\bigskip

\section{Identification of the limit law}
The identification of the limit  of any weakly converging sequence
of laws of local times of  fractional Brownian motions  will be a
consequence of some general results. The first one is inspired by
the occupation formula.

\begin{prop}\label{prop1idlim}
Let $\{X^n\}_{n\in\N}$ be a family of stochastic processes
verifying
\begin{enumerate}[(a)]

\item $\{X^n\}_{n\in\N}$  converges in law to $X$ in
$\mathcal{C}([0,T])$, when $n\to +\infty$.

\item Both the family
$\{X^n\}_{n\in\N}$ and $X$ have local times  $L^n$ and $L$
respectively, jointly continuous in $x$ and $t$.

\item The family of local times  $L^n$
converges in law to a process $Y$ in $\mathcal{C}([-D,D]\times
[0,T])$, when $n\to \infty$.
\end{enumerate}
Let $g:\R\times[-D,D]\to \R$ continuous with compact support such
that there exist $\alpha\in(0,1]$ and $C>0$ for which
\begin{equation}\label{condi}
\sup_{\substack{x\in[-D,D]\\
y,z\in\R}}\frac{|g(y,x)-g(z,x)|}{|y-z|^{\alpha}}<C.
\end{equation}
Then
\begin{equation*}
\int_{\R} g(u,x)Y(u,t)du\overset{\mathcal{L}}{=}\int_0^t g(X_s ,x)
ds,
\end{equation*}
in $\,\mathcal{C}([-D,D]\times [0,T])$.
\end{prop}
\begin{proof}
Fix $g$ continuous with compact support and satisfying condition
(\ref{condi}). Consider the maps
 $T_g$ and $U_g$, defined in the following way
\begin{align*}
T_g:\mathcal{C}(\R\times [0,T])&\longrightarrow
\mathcal{C}([-D,D]\times
[0,T])\\
y&\longmapsto T_g(y)(x,t)=\int_{\R} g(u,x)y(u,t)du,
\end{align*}
\begin{align*}
U_g:\mathcal{C}([0,T])&\longrightarrow \mathcal{C}([-D,D]\times
[0,T])\\
f&\longmapsto U_g (f)(x,t)=\int_0^t g(f(s),x)ds.
\end{align*}

It is easily checked that $T_g$ and $U_g$ are continuous maps with
respect to the usual topologies.

Proposition \ref{prop12} implies that, for any $(x,t)$ and
$n\in\N$,
\begin{equation}\label{igformocu}
\int_\R g(u,x)L_u^{t,n}du=\int_0^t g(X_s^n,x)ds.
\end{equation}
Due to the continuity of $\,T_g$ and $\,U_g$ and the convergence in
law of the families $\{X^n\}_{n\in \N}$ and $\{L^n\}_{n\in\N}$ to
$X$ and $Y$ in the spaces $\mathcal{C}([0,T])$ and
$\mathcal{C}([-D,D]\times[0,T])$, respectively, we  have
\begin{equation*}
\int_\R g(u,x)L_x^{t,n}du\stackrel{\law}{\longrightarrow}\int_\R
g(u,x)Y(u,t)du
\end{equation*}
and
\begin{equation*}
\int_0^t g(X_s^n,x)ds\stackrel{\law}{\longrightarrow}\int_0^t
g(X_s,x)ds,
\end{equation*}
in $\mathcal{C}([-D,D]\times [0,T])$, as $n\to \infty$.

Taking into account this convergence and using (\ref{igformocu})
we obtain
$$\int_\R g(u,x)Y(u,t)du\overset{\law}{=}\int_0^t
g(X_s,x)ds.$$ This concludes the proof.
\end{proof}

In the next proposition, we will prove that, under the hypotheses
of this last result, the finite dimensional distributions of $Y$
coincide with those of $L$.

\begin{prop}\label{propdisdimfin}
Let $\{X_n\}_{n\in\N}$ be a family of processes satisfying $(a)$,
$(b)$ and $(c)$ of the above proposition. Then, for any
$(x_1,t_1)$,..., $(x_k,t_k)\in[-D,D]\times [0,T]]$
\begin{equation*}
(Y(x_1,t_1),\dots,Y(x_k,t_k))\overset{\law}{=}(L_{x_1}^{t_1},\dots,
L_{x_k}^{t_k} ).
\end{equation*}
\end{prop}
\begin{proof}

Let $\varphi\in\mathcal C^1$ with compact support contained in
$[-1,1]$ and such that $\int_{\R}\varphi(x)dx=1$. Define, for any
 $\varepsilon>0$ and any $(u,x)\in \R\times [-D,D]
$,  $g_\varepsilon (u,x)=\frac{1}{\varepsilon}\varphi
\left(\frac{u-x}{\varepsilon}\right)$.

By Proposition \ref{prop1idlim}, we have that the following random
vectors are equal in law in $\R^k$
\begin{align*}
&\left(\int_\R g_\varepsilon(u,x_1)Y(t_1,u)du,\dots,\int_\R\!\!
g_\varepsilon(u,x_k)Y(t_k,u)du\right)\nonumber\\&\qquad\overset{\law}{=}\left(\int_0^{t_1}\!\!
g_\varepsilon(X_s,x_1)ds,\dots,\int_0^{t_k}\!\!g_\varepsilon(X_s,x_k)ds\right).
\end{align*}

Since $\{g_{\varepsilon}\}$  is an approximation of the identity,
for any fixed  $x\in\R$, $t\in [0,T]$ and any $\omega\in\Omega$ we
have that

\begin{equation*}
\lim_{\varepsilon \to 0}\int_\R g_\varepsilon(u,x)Y(u,t)du=Y(x,t),
\end{equation*}
and this implies that
\begin{equation*}
\left(\int_\R g_\varepsilon(u, x_1)Y(t_1,u)du,\dots,\int_\R\!\!
g_\varepsilon(u,
x_k)Y(t_k,u)du\right)\overset{\law}{\longrightarrow}\left(Y(x_1,t_1)\dots,Y(x_k,t_k)\right),
\end{equation*}
as $\varepsilon$ tends to $0$.

Using again that  ${g_{\varepsilon}}$ is an approximation of the
identity we obtain
\begin{equation*}
\left(\int_\R g_\varepsilon(u, x_1)L_u^{t_1}du,\dots,\int_\R\!\!
g_\varepsilon(u,
x_k)L_u^{t_k}du\right)\overset{\law}{\longrightarrow}\left(L_{x_1}^{t_1},\dots,L_{x_k}^{t_k}\right),
\end{equation*}
or equivalently, by using Proposition \ref{prop12}
\begin{equation*}
 \left(\int_0^{t_1}g_\varepsilon(X_s, x_1)ds, \dots
, \int_0^{t_k}g_\varepsilon(X_s,
x_k)ds\right)\overset{\law}{\longrightarrow}\left(L_{x_1}^{t_1},\dots,L_{x_k}^{t_k}\right),
\end{equation*}
when $\varepsilon\to 0$.

>From this, we conclude
\begin{equation*}
(Y(x_1,t_1),\dots,Y(x_k,t_k))\overset{\law}{=}(L_{x_1}^{t_1},\dots,L_{x_k}^{t_k}).
\end{equation*}

\end{proof}

The above general result can be applied to the fractional Brownian
motions.  Using also  the results of the preceding sections we
obtain the desired convergence in law of the family of local
times.

\begin{cor}
Given $H_0\in(0,1)$, the family $\{L^H\}_{H\in(0,1)}$ of local
times  of the fractional Brownian motions   converges in law to
the local time  $L^{H_0}$ of $B^{H_0}$ in
$\mathcal{C}([-D,D]\times[0,T])$, for any $D,\,T>0$ , when $H$
tends to $H_0$.
\end{cor}
\begin{proof}
By Proposition \ref{ajustprev} we have the tightness of the laws
of the family  $\{L^H\}_{H\in(H_0-\eta,H_0+\eta)}$, for certain
$\eta>0$, in $\mathcal{C}([-D,D]\times[0,T])$.

Take a sequence $\{H_n\}_{n\in \N}$ converging to $H_0$ as $n\to
\infty$ such that
\begin{equation}\label{conlawtl}
 L_x^{t,H_n}\stackrel{\mathcal{L}}{\longrightarrow} Y(x,t),
\end{equation}
in $\mathcal{C}([-D,D]\times [0,T])$, as $n\to \infty$.

Proposition \ref{propdisdimfin} gives that for any fixed
$(x_1,t_1)$, ..., $(x_k,t_k)$ we have that
\begin{equation*}
(Y(x_1,t_1),\dots,Y(x_k,t_k))\overset{\law}{=}(L_{x_1}^{H_0,\;t_1},\dots,L_{x_k}^{H_0,\;t_k}).
\end{equation*}
So, $\mathcal{L}(Y)=\mathcal{L}(L^{H_0})$ in
$\mathcal{C}([-D,D]\times[0,T])$.

\end{proof}

\section{Appendix}
\noindent{\it Proof of Theorem \ref{extlocal}}.

Fix $t\in [0,T]$. By condition $(ii)$, the  Fourier transform of
the occupation measure belongs to $L^2(\R)$, $\omega$-a.s.

By applying Theorem \ref{plancherel}, we have that $\omega$- a.s.
the measure $\mu^t$ is absolutely continuous and its density
$f^t\in L^2(\R)$ is square integrable in $x$. Moreover, defining
\begin{equation*}
f_N^t(x)=\frac{1}{2\pi}\int_{-N}^N e^{-iux}\phi^t(u)du=
\frac{1}{2\pi}\int_{-N}^N\int_0^t e^{-iux}e^{iu X_r}drdu,
\end{equation*}
we have
\begin{equation}\label{limdens11}
f_N^t\overset{L^2(\R)}{\longrightarrow} f^t,\qquad \omega-a.s.
\end{equation}

By Theorem \ref{thmabermanl2}, for any  $(x,t)$ there exists a
random variable, $L_x^t$, such that
\begin{equation*}
\sup_{(x,t)\in\R\times[0,T]}E|\psi_N(x,t)-L_x^t|^2\underset{N\to\infty}{\longrightarrow}0,
\end{equation*}
where
$\psi_N(x,t)=\frac{1}{2\pi}\int_{-N}^Ne^{-iux}\int_0^te^{iuX_r}drdu=f_N^t(x)$.

We will see that $\omega$-a.s. $L_{\cdot}^t$ is a version of the
density $f^t$ of $\mu^t$. Due to the above uniform convergence, in
any
  $[-A,A]$, with $A>0$, the following
convergence follows
\begin{equation*}
E\left(\int_{-A}^A|\psi_N(x,t)-L_x^t|^2dx\right)\underset{N\to\infty}\longrightarrow0.
\end{equation*}

That is,
\begin{equation*}
\int_{-A}^A|\psi_N(x,t)-L_x^t|^2dx\overset{L^1(\Omega)}{\longrightarrow}0,
\end{equation*}
when $N\to \infty$.  So, there exists a subsequence
$\{N_k\}_{k\in\N}$ and $\Omega'\subset \Omega$ with $P(\Omega')=1$
such that for any  $\omega\in \Omega'$, the integral
\begin{equation}\label{equax11}
\int_{-A}^A |\psi_{N_k}(x,t)-L_x^t|^2dx\longrightarrow 0,
\end{equation}
converges to 0, as $N_k\to \infty$.

On the other hand, by (\ref{limdens11}) we have, with probability
$1$,  that
\begin{equation}\label{equax12}
\int_{-\infty}^\infty |\psi_N(x,t)-f^t(x)|^2 dx\longrightarrow 0,
\end{equation}
when $N\to \infty$.

Finally, from (\ref{equax11}) and (\ref{equax12}), we can deduce
that  $L_x^t=f^t(x)$ $x$-a.e,  with probability $1$. So,
$\{L_x^t\}$ is a local time for
$X$.$\phantom{xxxxxxxxxxxxxxxxxxxxx}$ $\Box$

\section*{Acknowledgement}
This research was partially supported by DGES Grants BFM2003-01345
and BFM2003-00261

\bibliographystyle{empty}
\bibliography{}

\def\cprime{$'$}

\end{document}